\documentstyle[amssymb,amsfonts, 12pt]{amsart}

\newenvironment{proof}{\noindent {\bf Proof.}
}{\endprf\par}
\def \endprf{\hfill  {\vrule height6pt width6pt
depth0pt}\medskip}

\def\a{\alpha}
\def\b{\beta}
\def\e{\eta}
\def\ep{\varepsilon}
\def\Ga{\Gamma}
\def\la{\lambda}

\def\opsi{\overline{\psi}}
\def\p{\partial}

\begin{document}

 \theoremstyle{plain}
  \newtheorem{Theorem}{Theorem}
  \newtheorem{proposition}[subsection]{Proposition}
  \newtheorem{Lemma}{Lemma}
  \newtheorem{corollary}[subsection]{Corollary}

  \newtheorem{remark}[subsection]{Remark}
  \newtheorem{remarks}[subsection]{Remarks}

  \newtheorem{definition}[subsection]{Definition}

\title[Hypersurfaces with Prescribed $m-$th Mean Curvature in Hyperbolic Space]
{Starshaped compact hypersurfaces with prescribed
$m-$th Mean Curvature 
 in Hyperbolic Space}
\author{Qinian Jin and Yanyan Li
        }
\address{Rutgers University, New Brunswick, New Jersey
08903}
\email{yyli@@math.rutgers.edu, qjin@@math.rutgers.edu}

\maketitle

\baselineskip 17pt


\section{Introduction}

\noindent Let ${\Bbb S}^n$ be the unit sphere in the
Euclidean
space ${\Bbb R}^{n+1}$, and let $e$ be the standard
metric on
${\Bbb S}^n$ induced from ${\Bbb R}^{n+1}$. Suppose
that $(u,
\rho)$ are the spherical coordinates in ${\Bbb
R}^{n+1}$, where
$u\in {\Bbb S}^n$, $\rho\in [0, \infty)$. By choosing
the smooth
function $\varphi(\rho):=\sinh^2\rho$ on $[0,\infty)$
we can
define a Riemannian metric $h$ on the set $\{(u,
\rho): u\in {\Bbb
S}^n, 0\le \rho<\infty\}$ as follows
$$
h=d\rho^2+\varphi(\rho) e.
$$
This gives the space form ${\cal R}^{n+1}(-1)$ which
is the
hyperbolic space ${\Bbb H}^{n+1}$ with sectional
curvature $-1$.
For a smooth hypersurface ${\cal M}$ in ${\cal
R}^{n+1}(-1)$, we
denote by $\la_1, \cdots, \la_n$ its principal
curvatures with
respect to the metric $g:=h|_{\cal M}$. Then, for each
$1\le k\le
n$, the $k$-th mean curvature of ${\cal M}$ is defined
as
$$
H_k=\left(\begin{array}{lll}
n\\k\end{array}\right)^{-1}
\sum_{i_1<\cdots <i_k} \la_{i_1}\cdots\la_{i_k}.
$$

Let $\psi(u, \rho)$, $u\in {\Bbb S}^n$, $\rho\in
(0,\infty)$, be a
given positive smooth function satisfying suitable
conditions. We
are interested in the existence of a smooth
hypersurface ${\cal
M}$ embedded in ${\cal R}^{n+1}(-1)$ as a graph over
${\Bbb S}^n$
so that its $k$-th mean curvature is given by $\psi$.
We refer the
readers to \cite{CNS2} and
 \cite{BLO} for the introductory material
and the
history of this problem.

It is clear that ${\cal M}:=\{(u, z(u)): u\in {\Bbb
S}^n\}$ is an
embedded hypersurface in ${\cal R}^{n+1}(-1)$ for any
smooth
positive function $z$ on ${\Bbb S}^n$. We call $z$
$k$-admissible
if the principal curvatures $(\la_1(z(u)),\cdots,\la_n(z(u)))$ of
${\cal M}$
belong to $\Ga_k$, where $\Ga_k$ is the connected
component of
$\{\la\in {\Bbb R}^n: H_k(\la)>0\}$ containing the
positive cone
$\{\la\in {\Bbb R}^n: \la_1>0, \cdots, \la_n>0\}$.

The main result of this paper is the following

\begin{Theorem}\label{t1}
Let $1\le k\le n$, and let $\psi$ be a smooth positive
function in
the annulus $\overline{\Omega}: u\in {\Bbb S}^n,
\rho\in [R_1,
R_2], 0<R_1<R_2<\infty$, satisfying the conditions
$$
\psi(u, R_1)\ge \coth^k (R_1) \quad \mbox{and}\quad
\psi(u,
R_2)\le \coth^k (R_2)\quad \mbox{for~} u\in {\Bbb S}^n
$$
and
$$
\frac{\p}{\p \rho} \left(\psi(u,\rho) \sinh^k
\rho\right)\le 0
\quad \mbox{for all~} u\in {\Bbb S}^n \mbox{~and~}
\rho\in [R_1,
R_2].
$$
Then there exists a positive smooth $k$-admissible
function $z$ on
${\Bbb S}^n$ such that 
the closed
hypersurface ${\cal M}:=\{(u, z(u)): u\in {\Bbb
S}^n\}$ is in $\Omega\subset {\cal
R}^{n+1}(-1)$, and its 
 $k$-th mean curvature
is given by $\psi$:
$$
H_k(\la_1(z(u)),\cdots,\la_n(z(u)))=\psi(u, z(u))\quad
\forall\ u\in {\Bbb
S}^n.
$$
\end{Theorem}

In the Euclidean space (${\cal R}^{n+1}(0)$),  such results
were obtained in the case $k=0$ by Bakelman and Kantor
\cite{BK1}, \cite{BK2} and
by Treibergs and Wei \cite{TW}, in the case $k=n$ by 
Oliker \cite{O0}, and for general $k$ by Caffarelli,
Nirenberg and Spruck \cite{CNS2}.
In the elliptic space (${\cal R}^{n+1}(+1)$),  
such result is the combination of
 the work
of 
  Barbosa,  Lira
and Oliker
\cite{BLO} and that of  Li and Oliker \cite{LO}. 
The $k=n$ case in Theorem \ref{t1} was established
by Oliker  \cite{O}.
Our proof of 
Theorem \ref{t1}
uses the $C^0$ and $C^1$ a priori estimates obtained
in \cite{BLO} and  the arguments in \cite{LO} which is
based on the degree theory for fully nonlinear
elliptic operator of second order developed in \cite{L}. 
The main work for us to prove Theorem \ref{t1}
is to give the $C^2$ a priori estimates.
In establishing the $C^2$ estimates, we make use
Lemma \ref{L2},  a quantitative version of a theorem of 
Davis \cite{D} which, to our knowledge,
was given in 
\cite{A1}.  The theorem in \cite{D} says that a rotationally 
invariant function on symmetric matrices is concave
if and only if it is concave on the diagonal matrices,
while Lemma \ref{L2} allows the use of this term in
making $C^2$ a priori estimates.  The use of such a concave term
in $C^2$ estimates for solutions of
the Monge-Amp\`ere equation has been extensive,
see e.g. Calabi \cite{Ca1} and
Pogorelov \cite{P}.
The use of Lemma  \ref{L2} in $C^2$ estimates for solutions of
more general equations can be found in \cite{A1},
\cite{A2}, \cite{G}, \cite{SUW}, \cite{U1} and
\cite{U2}.

\bigskip

\noindent{\bf Acknowledgement.}\  The work of the second author is 
partially supported by
       NSF grant DMS-0401118.

\section{Some fundamental formulae}

\noindent
Let us define a function $f$ on $\Ga_k$ by
$$
f(\la)=\left[\left(\begin{array}{lll}
n\\k\end{array}\right)^{-1}
\sum_{i_1<\cdots <i_k}
\la_{i_1}\cdots\la_{i_k}\right]^{1/k},
$$
where $\la:=(\la_1,\cdots, \la_n)\in \Ga_k$. It is
well known that
$f$ is smooth, positive, concave, and strictly
increasing with
respect to each variable, see e.g. \cite{CNS}. Now our problem is
equivalent to finding
a smooth positive $k$-admissible function $z$ on
${\Bbb S}^n$ so
that
\begin{equation}\label{1}
F({\bf B})=\opsi
\end{equation}
on ${\cal M}:=\{(u,z(u)): u\in {\Bbb S}^n\}$, where
$\opsi:=\psi^{1/k}$, ${\bf B}$ is the second
fundamental form of
${\cal M}$, and $F({\bf B}):=f(\la({\bf B}))$ with
$\la({\bf B})$
being the eigenvalues of ${\bf B}$ with respect to the
metric $g$
on ${\cal M}$.

Suppose now ${\cal M}$ is the graph of a smooth
positive
$k$-admissible function $z$ on ${\Bbb S}^n$. Let us
recall the
formulae given in \cite{BLO} for the components of
$g=(g_{ij})$
and ${\bf B}=(b_{ij})$ on ${\cal M}$ under a local
coordinate. Let
$\theta^1, \cdots, \theta^n$ be a smooth local
coordinate of
${\Bbb S}^n$, which of course gives a local coordinate
of ${\cal
M}$. If we denote by $\{e_{ij}\}$ the components of
$e$ under this
local coordinate, and set $z_i=\frac{\p z}{\p
\theta^i}$ and
$z_{ij}=\frac{\p^2z}{\p\theta^i\p\theta^j}$, then
\begin{equation}\label{2}
g_{ij}=\varphi e_{ij}+z_iz_j
\end{equation}
and
\begin{equation}\label{3}
g^{ij}=\frac{1}{\varphi}\left[e^{ij}-\frac{z^iz^j}{\varphi+|\nabla'
z|^2}\right], \qquad z^i=e^{ij} z_j,
\end{equation}
where $(g^{ij})=(g_{ij})^{-1}$,
$(e^{ij})=(e_{ij})^{-1}$, and
$\nabla'$ denotes the Levi-Civita connection on ${\Bbb
S}^n$.
Moreover, for the second fundamental form we have
\begin{equation}\label{3.1}
b_{ij}=\frac{\varphi}{\sqrt{\varphi^2+\varphi|\nabla'
z|^2}}\left\{-\nabla'_{ij} z+\frac{\p \ln \varphi}{\p
\rho}z_iz_j+\frac{1}{2} \frac{\p \varphi}{\p \rho}
e_{ij}\right\}.
\end{equation}

We also need the following well-known fundamental
equations for a
hypersurface ${\cal M}$ in ${\cal R}^{n+1}(-1)$:
\begin{eqnarray}
\mbox{Codazzi equation:} &&\nabla_i b_{jk}=\nabla_j
b_{ki}=\nabla_k b_{ij}\label{3.2}\\
\mbox{Gauss equation:} &&
R_{ijkl}=(b_{ik}b_{jl}-b_{il}b_{jk})-(g_{ik}g_{jl}-g_{il}g_{jk})
\label{3.3}\\
\mbox{Ricci equation:} && \nabla_l\nabla_k
b_{ij}-\nabla_k\nabla_l
b_{ij}= b_{ip}g^{pq} R_{qjkl} + b_{jp}g^{pq} R_{qikl}
\label{3.4}
\end{eqnarray}
where $R_{ijkl}$ denotes the Riemannian curvature
tensor of ${\cal
M}$, and $\nabla_i$ and $\nabla_i\nabla_j$ the
covariant
differentiations in the metric $g$ with respect to
some local
coordinates on ${\cal M}$.

As the preparation for deriving the $C^2$-estimates,
let us
introduce the following two functions on ${\cal M}$
\begin{equation}\label{3.5}
\tau=
\frac{\varphi(z)}{\sqrt{\varphi(z)+|\nabla'z|^2}}
\quad
\mbox{and}\quad \e=-\cosh(z).
\end{equation}
We have

\begin{Lemma}\label{L1}
For $\tau$ and $\e$ the following equations hold
\begin{eqnarray}
\nabla_i \tau &=& -b_{ip}g^{pq} \nabla_q
\e,\label{3.6}\\
\nabla_{ij}\tau &=& -\nabla_p b_{ij} g^{pq} \nabla_q
\e
-\tau b_{ip} g^{pq} b_{qj}-\e b_{ij}, \label{3.7}\\
\nabla_{ij} \e &=& \tau b_{ij} +\e g_{ij}.\label{3.8}
\end{eqnarray}
\end{Lemma}
\begin{proof}  These formulae have been derived in
\cite{BLO} by using
another model of ${\cal R}^{n+1}(-1)$. In fact, we can
show them
directly. Since (\ref{3.7}) is an immediate
consequence of
(\ref{3.6}), (\ref{3.8}) and the Codazzi equation
(\ref{3.2}), it
suffices to verify (\ref{3.6}) and (\ref{3.8}). Let
$c(\rho)=\cosh(\rho)$ and $s(\rho)=\sinh(\rho)$. Then
\begin{eqnarray*}
\nabla_i\tau &=& \frac{2sc}{\sqrt{\varphi+|\nabla'
z|^2}}z_i-\frac{\varphi}{(\varphi+|\nabla'z|^2)^{3/2}}\left(sc
z_i+e^{pq}\nabla_{ip}'z \nabla_q' z\right)\\
&=& \frac{1}{(\varphi+|\nabla'z|^2)^{3/2}}\left\{sc
\varphi
z_i+2scz_i|\nabla'z|^2-\varphi \nabla_{ip}'z
z^p\right\}
\end{eqnarray*}
Noting that $\nabla_q \e=- s z_q$, we have from
(\ref{3}) and
(\ref{3.1}) that
$$
b_{ip}g^{pq}
\nabla_q\e=-\frac{s}{\varphi+|\nabla'z|^2} b_{ip}
z^p=-\nabla_i\tau.
$$

Let us now verify (\ref{3.8}) for any fixed
$\bar{u}\in {\Bbb
S}^n$. Noting that the both sides of (\ref{3.8}) are
tensorial, we
may assume that the local coordinates are chosen such
that
$\frac{\p g_{jk}}{\p \theta^i}=0$ at $\bar{u}$. Then
from
(\ref{2}) we have
$$
\frac{\p e_{lj}}{\p \theta^i}
=-\frac{2sc}{\varphi^2}\left(g_{lj}-z_lz_j\right)
z_i-\frac{1}{\varphi}
\left(z_{il}z_j+z_lz_{ij}\right).
$$
Thus the corresponding Christoffel symbols of ${\Bbb
S}^n$ are
given by
\begin{eqnarray}\label{3.9}
{\Gamma'}_{ij}^k&=&\frac{1}{2} e^{kl}\left\{\frac{\p
e_{lj}}{\p
\theta^i}+\frac{\p e_{li}}{\p \theta^j}-\frac{\p
e_{ij}}{\p
\theta^l}\right\}\nonumber\\
&=& -\frac{sc}{\varphi^2}\left\{e^{kl}g_{lj}z_i +
e^{kl}g_{li}z_j
-g_{ij}z^{k}-z_i z_jz^k\right\}
-\frac{1}{\varphi}z_{ij}z^k.
\end{eqnarray}
This, together with (\ref{2}), gives
\begin{eqnarray*}
\nabla'_{ij} z&=& z_{ij}-{\Gamma'}_{ij}^k z_k\\
&=& \frac{\varphi+|\nabla'z|^2}{\varphi} z_{ij} +
\frac{sc}{\varphi^2} \left\{g_{lj}z_iz^l +g_{li}z_j
z^l-
g_{ij}|\nabla'z|^2-z_iz_j
|\nabla'z|^2\right\}\\
&=& \frac{\varphi+|\nabla'z|^2}{\varphi} z_{ij}
+\frac{sc}
{\varphi} \left(2z_iz_j-|\nabla'z|^2 e_{ij}\right).
\end{eqnarray*}
Noting that $\nabla_{ij}z=z_{ij}$ at $\bar{u}$, we
thus have
\begin{equation}\label{3.10}
\nabla_{ij}z=\frac{1}{\varphi+|\nabla'z|^2}\left(\varphi\nabla'_{ij}
z-2sc z_iz_j +sc |\nabla' z|^2 e_{ij}\right).
\end{equation}
Therefore
$$
\nabla_{ij}\e = -cz_iz_j-s\nabla_{ij}
z=-cz_iz_j-\frac{s}{\varphi+|\nabla'z|^2}\left(\varphi\nabla'_{ij}
z-2sc z_iz_j+sc |\nabla' z|^2 e_{ij}\right).
$$
But from (\ref{2}) and (\ref{3.1}) we can see that the
right hand
side of the above equation is exactly $\tau b_{ij}+\e
g_{ij}$.
\end{proof}

\section{$C^2$-estimates}

\noindent Now we are in a position to derive the $C^2$
estimates
for any smooth positive $k$-admissible solutions of
(\ref{1}) in
${\cal R}^{n+1}(-1)$. Let us set
$$
f_i=\frac{\p f}{\p\lambda_i}, \quad  F^{ij}=\frac{\p
F}{\p b_{ij}}
\quad \mbox{and} \quad F^{ij, kl}=\frac{\p^2F}{\p
b_{ij}\p
b_{kl}}.
$$
We will achieve our aim by choosing suitable test
function and
making full use of the terms involving $F^{ij, kl}$.
In
particular, we need the following

\begin{Lemma}
(\cite{A1})\label{L2}  \
For any symmetric matrix $(\e_{ij})$ there holds
$$
F^{ij, kl}\e_{ij}\e_{kl}=\sum_{i,j} \frac{\p^2 f}{\p
\la_i \p
\la_j} \e_{ii} \e_{jj} +\sum_{i\ne j}
\frac{f_i-f_j}{\la_i-\la_j}
\e_{ij}^2,
$$
where the second term on the right-hand side must be
interpreted
as a limit whenever $\la_i=\la_j$.
\end{Lemma}

This result was, to our knowledge,
first stated in \cite{A1};
for proofs
one may consult \cite{G, A2}.

\begin{Theorem}\label{t2}
Let $1\le k\le n$ and let $\psi$ be a positive $C^2$
function in
the annulus $\overline{\Omega}: u\in {\Bbb S}^n,
\rho\in [R_1,
R_2], 0<R_1<R_2<a$. Let $z\in C^4({\Bbb S}^n)$ be a
positive
$k$-admissible solution of (\ref{1}) in ${\cal
R}^{n+1}(-1)$
satisfying
$$
R_1\le z\le R_2 \quad \mbox{and} \quad |\nabla' z|\le
C_0=constant
\quad \mbox{on~~} {\Bbb S}^n.
$$
Then
$$
\|z\|_{C^2({\Bbb S}^n)}\le C,
$$
where the constant $C$ depends only on $k$, $n$,
$R_1$, $R_2$,
$C_0$ and $\|\psi\|_{C^2(\overline{\Omega})}$.
\end{Theorem}
\begin{proof} We will estimate the maximal principal
curvature of
${\cal M}$. Since $z$ is $k$-admissible, this
estimate, together
with the $C^0$ and $C^1$ bounds of $z$ and the
equation
(\ref{3.1}), implies an estimate for $\|z\|_{C^2({\Bbb
S}^n)}$.
Consider the function
$$
\widetilde{W}(u, \xi)={\bf B}(\xi,\xi)
\exp\left[\Phi(\tau)-\beta
\eta\right],
$$
where $u\in {\Bbb S}^n$, $\xi$ is a unit tangent
vector of ${\cal
M}$ at $(u, z(u))$, $\tau$ and $\eta$ are defined as
in
(\ref{3.5}), and the function $\Phi$ and the constant
$\b>0$ will
be determined later. Suppose the maximum of
$\widetilde{W}$ is
attained at some point $\bar{u}\in {\Bbb S}^n$ in the
unit
tangential direction $\bar{\xi}$ of ${\cal M}$ at
$(\bar{u},
z(\bar{u}))$. We may choose the local coordinates
$\theta^1,
\cdots, \theta^n$ around $\bar{u}$ such that
$$
g_{ij}=\delta_{ij}\,\,\, \mbox{~~~and~~~}\,\,\,
\frac{\p
g_{ij}}{\p \theta^k}=0 \mbox{~~at~} \bar{u}.
$$
Moreover, since $\bar{\xi}$ is the maximal principal
direction of
${\cal M}$ at $(\bar{u}, z(\bar{u}))$, such
coordinates can be
chosen so that $\{b_{ij}\}$ is diagonal at $\bar{u}$
and
$b_{11}(\bar{u})={\bf B}(\bar{\xi},\bar{\xi})$.

Consider the local function $Z=b_{11}/g_{11}$. By
direct
calculation we have at $\bar{u}$ that
$$
\nabla_iZ=\frac{\p b_{11}}{\p \theta^i}=\nabla_i
b_{11}
$$
and
$$
\nabla_i\nabla_j Z=\frac{\p^2 b_{11}}{\p \theta^i\p
\theta^j}-b_{11} \frac{\p^2 g_{11}}{\p \theta^i\p
\theta^j}
=\frac{\p^2 b_{11}}{\p \theta^i\p \theta^j}-2\frac{\p
\Gamma_{j1}^1}{\p \theta^i} b_{11}=\nabla_i\nabla_j
b_{11}.
$$
It is clear that the function
$$
W(u)=Z(u)\exp\left[\Phi(\tau)-\beta\eta\right].
$$
has a local maximum at $\bar{u}$. Thus at $\bar{u}$
\begin{equation}\label{5}
0=\nabla_i(\log W)=\frac{\nabla_i Z}{Z}+\Phi'\nabla_i
\tau -\b
\nabla_i \e=\frac{\nabla_i
b_{11}}{b_{11}}+\Phi'\nabla_i \tau -\b
\nabla_i \e
\end{equation}
and the matrix
\begin{eqnarray*}
\left\{\nabla_{ij}(\log
W)\right\}&=&\left\{\frac{\nabla_i\nabla_j
Z }{Z}-\frac{\nabla_i Z\nabla_j Z}{Z^2} +\Phi'
\nabla_{ij} \tau
+\Phi'' \nabla_i\tau \nabla_j\tau -\b \nabla_{ij}
\e\right\}\\
&=&\left\{\frac{\nabla_i\nabla_j
b_{11}}{b_{11}}-\frac{\nabla_i
b_{11}\nabla_j b_{11}}{b_{11}^2} +\Phi' \nabla_{ij}
\tau +\Phi''
\nabla_i\tau \nabla_j\tau -\b \nabla_{ij} \e\right\}
\end{eqnarray*}
is negative semi-definite. Therefore
\begin{align}\label{6}
0\ge F^{ij}\nabla_{ij}(\log W)=& \frac{1}{b_{11}}
F^{ij}
\nabla_i\nabla_j b_{11}-\frac{1}{b_{11}^2} F^{ij}
\nabla_i
b_{11}\nabla_j b_{11} +\Phi'
F^{ij} \nabla_{ij} \tau \nonumber\\
& +\Phi'' F^{ij} \nabla_i\tau\nabla_j \tau -\beta
F^{ij}\nabla_{ij} \eta.
\end{align}
Since $\{b_{ij}\}$ is diagonal at $\bar{u}$,
$\{F^{ij}\}$ is also
diagonal there and $F^{ii}=f_i$. For simplicity, we
let $\la_i
=b_{ii}(\bar{u})$ and assume $\la_1\ge \la_2\ge \cdots
\ge \la_n$,
moreover we may assume $\la_1\ge 1$. Then, see lemma 2
in \cite{EH} or lemma A.2 in \cite{L2},
we have $f_1\le f_2\le \cdots \le f_n$. It follows
from (\ref{6})
that
\begin{eqnarray}\label{7}
0&\ge& \frac{1}{\la_1}\sum_i f_i \nabla_i\nabla_i
b_{11}-\frac{1}{\la_1^2}\sum_i f_i |\nabla_i b_{11}|^2
+\Phi'
\sum_i f_i \nabla_{ii} \tau\nonumber\\
&&+ \Phi''\sum_i f_i |\nabla_i \tau|^2-\beta \sum_i
f_i
\nabla_{ii} \eta.
\end{eqnarray}
Now we take the covariant differentiation on (\ref{1})
to get
\begin{equation}\label{8}
F^{ij}\nabla_1\nabla_1 b_{ij} + F^{ij, kl} \nabla_1
b_{ij}\nabla_1
b_{kl} =\nabla_{11}\opsi.
\end{equation}
>From (\ref{3.2}), (\ref{3.3}) and (\ref{3.4}) it
follows that
\begin{eqnarray*}
\nabla_1\nabla_1 b_{ii}&=&\nabla_1\nabla_i
b_{1i}=\nabla_i\nabla_1
b_{1i}+\sum_k b_{1k} R_{kii1}+\sum_k b_{ik}
R_{k1i1}\nonumber\\
&=& \nabla_i\nabla_i b_{11}+b_{11}b_{ii}^2
-b_{11}^2b_{ii}
-(b_{11}\delta_{1i}
-b_{11}\delta_{ii}+b_{ii}-b_{ii}\delta_{1i}).
\end{eqnarray*}
This shows that
$$
F^{ij}\nabla_1\nabla_1 b_{ij}=\sum_i f_i
\nabla_i\nabla_i b_{11}
+\la_1\sum_i f_i \la_i^2 -\la_1^2 \sum_i f_i \la_i
+\la_1 {\cal
T}-\sum_i f_i \la_i,
$$
where ${\cal T}:=\sum_i f_i$. Since the degree one
homogeneity of
$f$ implies $\sum_i f_i\la_i=\opsi$, the above
equation together
with (\ref{8}) gives
$$
\sum_i f_i \nabla_i\nabla_i b_{11}=-F^{ij, kl}
\nabla_1
b_{ij}\nabla_1 b_{kl}+\nabla_{11}\opsi +\la_1^2 \opsi+
\opsi
-\la_1\sum_i f_i \la_i^2 -\la_1 {\cal T}.
$$
Plugging this into (\ref{7}), noting that $\opsi\ge
c_0>0$ and
$\la_1\ge 1$,  we therefore obtain
\begin{eqnarray}\label{9}
0&\ge& c_0\la_1-\frac{1}{\la_1} F^{ij, kl} \nabla_1
b_{ij}\nabla_1
b_{kl} +\frac{\nabla_{11}\opsi}{\la_1} -\sum_i f_i
\la_i^2 -{\cal
T} -\frac{1}{\la_1^2} \sum_i f_i |\nabla_i
b_{11}|^2 \nonumber\\
&& +\Phi' \sum_i f_i \nabla_{ii}\tau +\Phi''\sum_i f_i
|\nabla_i
\tau|^2 -\beta \sum_i f_i \nabla_{ii}\eta.
\end{eqnarray}
>From (\ref{3.7}) and (\ref{3.8}) we have
\begin{equation}\label{11}
\b \sum_i f_i\nabla_{ii}\e =\b \tau \sum_i f_i
\la_i+\b \e {\cal
T} =\b \tau \opsi +\b  \e {\cal T}.
\end{equation}
and
\begin{eqnarray}\label{12}
\Phi' \sum_i f_i \nabla_{ii}\tau &=&
\Phi'\left\{-\sum_p \nabla_p
\e \left(\sum_i f_i \nabla_p b_{ii}\right)-\tau \sum_i
f_i \la_i^2 - \e \sum_i f_i \la_i \right\} \nonumber\\
&=& \Phi'\left\{-\sum_p \nabla_p\e \nabla_p \opsi
-\tau \sum_i
f_i \la_i^2 -\e \opsi\right\}\nonumber\\
&\ge & -C|\Phi'|-\Phi' \tau \sum_i f_i \la_i^2.
\end{eqnarray}
Here we used the facts $|\nabla_p\e|\le C$ and
$|\nabla_p\opsi|\le
C$ at $\bar{u}$ which can be demonstrated as follows.
Since
$g_{ij}=\delta_{ij}$ at $\bar{u}$, it follows from
(\ref{2}) that
$(z_p)^2\le 1$ at $\bar{u}$. Note that $\nabla_p \e
=-\sinh(z)
z_p$. Therefore $|\nabla_i \e|\le C$ at $\bar{u}$. For
$|\nabla_p
\opsi|$, we note that $\nabla_p \opsi= \opsi_p+\opsi_z
 z_p$.
Thus, by using (\ref{3}), we have at $\bar{u}$ that
\begin{eqnarray*}
|\nabla_p \opsi|^2 &\le&
C\left(1+|\opsi_p|^2\right)\le
C\left(1+g^{ij} \opsi_i\opsi_j\right)\le
C\left(1+\varphi^{-1}
e^{ij}\opsi_i\opsi_j\right)\\
&=& C\left (1+\varphi^{-1}|\nabla'\opsi|^2\right)\le
C.
\end{eqnarray*}

One can show that
\begin{equation}\label{12.5}
\frac{\nabla_{11}\opsi}{\la_1}\ge -C \quad \mbox{at~}
\bar{u}.
\end{equation}
To see this, note that $\frac{\p g_{ij}}{\p
\theta^k}=0$ at
$\bar{u}$, we have
$$
\nabla_{11}\opsi=\opsi_{11}+2\opsi_{z1} z_1+\opsi_{zz}
\left(z_1\right)^2+\opsi_z z_{11}.
$$
Similar to the above argument we can show
$|\opsi_{z1}|\le C$ at
$\bar{u}$. Therefore at $\bar{u}$
$$
|\nabla_{11}\opsi|\le
C\left(1+|\opsi_{11}|+|z_{11}|\right).
$$
Let us estimate $|\opsi_{11}|$. It follows from
(\ref{3}) that
\begin{eqnarray*}
|\nabla'_{11}\opsi|^2&\le&
g^{ik}g^{jl}\nabla'_{ij}\opsi\nabla'_{kl}\opsi\\
&\le&\varphi^{-2}e^{ik}e^{jl}\nabla'_{ij}\opsi\nabla'_{kl}\opsi
+\varphi^{-2}(\varphi+|\nabla'z|^2)^{-2}\left(z^iz^j\nabla'_{ij}\opsi\right)^2\\
&\le& \varphi^{-2}|{\nabla'}^2\opsi|^2
+\varphi^{-2}(\varphi^2
+|\nabla'z|^2)^{-2}
|{\nabla'}^2\opsi|^2|\nabla'z \otimes \nabla'z|^2\\
&\le& \varphi^{-2}|{\nabla'}^2\opsi|^2
+\varphi^{-2}(\varphi^2
+|\nabla'z|^2)^{-2}
|{\nabla'}^2\opsi|^2 |\nabla'z|^4\\
&\le& C
\end{eqnarray*}
By using (\ref{2}) we obtain at $\bar{u}$ that
$$
\sum_i |z^i|^2=g_{ij}z^iz^j=\varphi
|\nabla'z|^2+|\nabla'z|^4\le
C.
$$
This together with (\ref{3}) then implies $|e^{ij}|\le
C$. Thus it
follows from (\ref{3.9}) that $|{\Gamma'}_{11}^k|\le
C(1+|z_{11}|)$ at $\bar{u}$. Since
$\opsi_{11}=\nabla'_{11}\opsi
+{\Gamma'}_{11}^k\opsi_k$, we therefore have
$|\opsi_{11}|\le C
\left(1+|z_{11}|\right)$. Since $z_{11}=\nabla_{11}z$
at
$\bar{u}$, from (\ref{3.10}) and (\ref{3.1}) we
finally obtain
$$
|\nabla_{11}\opsi|\le C\left(1+|z_{11}|\right)\le
C\left(1+\la_1\right)
$$
which gives (\ref{12.5}).

Combining (\ref{9}), (\ref{11}), (\ref{12}) and
(\ref{12.5}), we
thus obtain
\begin{align}\label{13}
0\ge&  c_0 \la_1-C(1+|\Phi'|)-(1+\b \e) {\cal T}- \b
\tau \opsi-
(\Phi' \tau +1)\sum_i f_i\la_i^2+\Phi'' \sum_i f_i
|\nabla_i\tau|^2\nonumber\\
& -\frac{1}{\la_1^2} \sum_i f_i |\nabla_i b_{11}|^2
-\frac{1}{\la_1}  F^{ij, kl} \nabla_1 b_{ij} \nabla_1
b_{kl}.
\end{align}
Now we will use Lemma \ref{L2}, similar to the way
used in  \cite{SUW}.

{\it Case 1.} $\lambda_n <-\theta \lambda_1$ for some
positive
constant $\theta$ (to be chosen later).

In this case, using the concavity of $F$ we may
discard the last
term on the right hand side of (\ref{13}) since it is
nonnegative.
Also from (\ref{5}) we have for any $\ep>0$
\begin{eqnarray*}
\frac{1}{\la_1^2} \sum_i f_i |\nabla_i b_{11}|^2 &=&
\sum_i f_i
|\Phi' \nabla_i \tau-\b \nabla_i \e|^2\\
&\le & (1+\ep^{-1}) \b^2 \sum_i f_i |\nabla_i \e|^2
+(1+\ep)
(\Phi')^2\sum_i f_i |\nabla_i \tau|^2.
\end{eqnarray*}
Therefore, from (\ref{13}) it yields
\begin{align}\label{14}
0&\ge c_0 \la_1 -C(1+|\Phi'|) -\left[(1+\b\e)
+C(1+\ep^{-1})\b^2\right] {\cal T} -\b \tau \opsi-
(\Phi'\tau+1)\sum_i f_i \la_i^2 \nonumber\\
& +\left[ \Phi''-(1+\ep)(\Phi')^2\right] \sum_i f_i
|\nabla_i
\tau|^2.
\end{align}
Using (\ref{3.6}) we have
$$
\sum_i f_i |\nabla_i \tau|^2=\sum_i f_i \la_i^2
|\nabla_i\e|^2 \le
c_1 \sum_i f_i \la_i^2
$$
for some constant $c_1>0$. If we can choose $\Phi$
such that
$\Phi''-(1+\ep) (\Phi')^2\le 0$, then from (\ref{14})
we have
\begin{eqnarray}\label{15}
0&\ge& c_0 \la_1- C(1+|\Phi'|)-\left[(1+\b
\e) +C(1+\ep^{-1})\b^2\right] {\cal T}-\b \tau
\opsi\nonumber\\
&& +\left[-(\Phi'\tau+1) +c_1 \left(\Phi''
-(1+\ep)(\Phi')^2\right)\right]\sum_i f_i\la_i^2.
\end{eqnarray}
In order to choose $\Phi$, let $a>0$ be a positive
number such
that $\tau \ge 2a$ which is guaranteed by our
assumption. Then we
define
$$
\Phi(\tau)=-\log(\tau-a).
$$
It is easy to check that $\Phi''-(1+\ep)(\Phi')^2<0$.
Moreover,
for $\ep=\frac{a^2}{2c_1}$ we have
$$
-(\Phi'\tau+1)+c_1\left(\Phi''-(1+\ep)(\Phi')^2\right)
=\frac{a}{\tau-a}-\frac{c_1\ep} {(\tau-a)^2}\ge
\frac{a^2}{2(\tau-a)^2} \ge c_2>0.
$$
Therefore we get from (\ref{15}) that
$$
0\ge c_0 \la_1-C-C{\cal T}+ c_2 \sum_i f_i \la_i^2.
$$
Since $\lambda_n \le -\theta \lambda_1$ and $f_n\ge
\frac{1}{n}{\cal T}$, we have $\sum_i f_i \la_i^2 \ge
f_n
\la_n^2\ge \frac{1}{n}\theta^2 {\cal T} \la_1^2$.
Hence
$$
0\ge c_0 \la_1-C-C{\cal T}+\frac{c_2\theta^2}{n}{\cal
T} \la_1^2.
$$
This clearly implies $\la_1$ is bounded from above.

{\it Case 2.} $\lambda_n\ge -\theta \lambda_1$.

We now have $\la_i\ge -\theta \la_1$ for all $i=1,
\cdots, n$. Let
us partition $\{1, \cdots, n\}$ into two parts:
$I=\{j: f_j\le
4f_1\}$ and $J=\{j: f_j>4f_1\}$. Using (\ref{5}) we
have for $i\in
I$ that
\begin{eqnarray*}
\frac{1}{\la_1^2} f_i |\nabla_i b_{11}|^2 &=& f_i
|\Phi' \nabla_i \tau-\b \nabla_i \e|^2\\
&\le & (1+\ep) (\Phi')^2 f_i |\nabla_i \tau|^2
+(1+\ep^{-1})\b^2 f_i |\nabla_i \e|^2\\
&\le& (1+\ep) (\Phi')^2 f_i |\nabla_i \tau|^2
+C(1+\ep^{-1}) \b^2
f_1.
\end{eqnarray*}
Therefore it follows from (\ref{13}) that
\begin{eqnarray*}
0&\ge& c_0 \la_1 -C(1+|\Phi'|)-\b\tau \opsi -(1+\b
\e){\cal T}
-(\Phi' \tau +1) \sum_i f_i \la_i^2\\
&& +\left[\Phi''-(1+\ep)(\Phi')^2 \right] \sum_i f_i
|\nabla_i\tau|^2 -C(1+\ep^{-1}) \b^2 f_1\\
&& -\frac{1}{\la_1^2} \sum_{j\in J} f_j|\nabla_j
b_{11}|^2
-\frac{1}{\la_1} F^{ij, kl} \nabla_1 b_{ij} \nabla_1
b_{kl}.
\end{eqnarray*}
Proceeding exactly as before we have
$$
-(\Phi'\tau+1) \sum_i f_i \la_i^2
+\left[\Phi''-(1+\ep)(\Phi')^2\right] \sum_i f_i
|\nabla_i \tau|^2
\ge c_2 \sum_i  f_i \la_i^2
$$
if we choose $\ep=\frac{a^2}{2c_1}$. So
\begin{eqnarray}\label{16}
0&\ge& c_0 \la_1-C(1+|\Phi'|)-\b \tau \opsi-(1+\b
\e){\cal T}
+c_2\sum_i  f_i \la_i^2-C(1+\ep^{-1})\b^2 f_1
\nonumber\\
&& -\frac{1}{\la_1^2} \sum_{j\in J} f_j |\nabla_j
b_{11}|^2
-\frac{1}{\la_1}F^{ij,kl}\nabla_1 b_{ij}\nabla_1
b_{kl}.
\end{eqnarray}
By using Lemma \ref{L2} and noting $1\not\in J$ we
have
$$
-\frac{1}{\la_1} F^{ij,kl} \nabla_1 b_{ij}\nabla_1
b_{kl} \ge
-\frac{2}{\la_1} \sum_{j\in J} \frac{f_1-f_j}
{\la_1-\la_j}
|\nabla_1 b_{1j}|^2 = -\frac{2}{\la_1} \sum_{j\in J}
\frac{f_1-f_j}{\la_1-\la_j} |\nabla_j b_{11}|^2.
$$
Therefore
\begin{eqnarray}\label{17}
0&\ge& c_0 \la_1 -C(1+|\Phi'|)-\b\tau\opsi-(1+\b\e)
{\cal
T} -C(1+\ep^{-1}) \b^2 f_1+c_2 \sum_i
f_i\la_i^2\nonumber\\
&& -\frac{2}{\la_1}\sum_{j\in J} \frac{f_1-f_j}
{\la_1-\la_j}
|\nabla_j b_{11}|^2 -\frac{1}{\la_1^2} \sum_{j\in J}
f_j |\nabla_j
b_{11}|^2.
\end{eqnarray}
We claim that
$$
-\frac{2}{\la_1} \frac{f_1-f_j}{\la_1-\la_j}\ge
\frac{1}{\la_1^2}
f_j, \qquad \forall j\in J.
$$
This is equivalent to showing $2f_1\la_1\le f_j
\la_1++f_j \la_j$.
Since $j\in J$, we have $f_j>4f_1$. If $\la_j \ge 0$,
this is
obviously true. If $\la_j<0$, then $-\theta\la_1\le
\la_j <0$, and
hence
$$
f_j\la_1+f_j\la_j\ge (1-\theta) f_j\la_1\ge
4(1-\theta)f_1\la_1\ge
2f_1 \la_1
$$
if we choose $\theta=\frac{1}{2}$. From this claim and
(\ref{17})
we obtain
$$
0\ge c_0 \la_1-C-\b\tau \opsi-(1+\b\e) {\cal T}
+c_2\sum_i  f_i
\la_i^2-C(1+\ep^{-1})\b^2 f_1.
$$
Recall the definition of $\e$, we have $-c_3\le
\eta\le -c_4$ for
two positive constants $c_3$ and $c_4$. Choose $\b$ to
be
sufficiently large so that $-(1+\b \e)\ge 0$. Then we
get
$$
0\ge -C+c_0 \la_1+c_2 f_1 \la_1^2- C f_1.
$$
This clearly implies an upper bound for $\la_1$.
\end{proof}

\section{Proof of main result}

\noindent  Since the proof of Theorem \ref{t1}
essentially follows
the lines in \cite{LO}, only the sketch will be given
below.

We may assume that neither $z(u)\equiv R_1$ nor
$z(u)\equiv R_2$
is a solution of (\ref{1}); otherwise we are done. Let
us fix some
$\overline{R}$ such that $R_1<\overline{R}<R_2$ and
define a
family of functions
$$
\opsi^t(u,\rho):=t\opsi(u,\rho)+(1-t)A^\ep
\coth^{1+\ep}(\rho),
\qquad t\in [0,1],
$$
where  $\ep>0$ is a positive constant and
$A=\coth^{-1}
(\overline{R})$. Fix $0<\a<1$, and denote by
$C_a^{4,\a}({\Bbb
S}^n)$ the subset of functions from $C^{4,\a}({\Bbb
S}^n)$ which
is $k$-admissible. We define a family of operators
$\Psi(\cdot,
t): C_a^{4,\a}({\Bbb S}^n)\to C^{2,\a}({\Bbb S}^n)$ by
$$
\Psi(z(u),t)\equiv F({\bf B})-\opsi^t(u, z(u)), \quad
u\in {\Bbb
S}^n,
$$
where $z\in C_a^{4,\a}({\Bbb S}^n)$ and ${\bf B}$ is
the second
fundamental form of ${\cal M}:=\{(u,z(u)): u\in {\Bbb
S}^n\}$.
Consider the family of equations
\begin{equation}\label{27}
\Psi(z,t)\equiv 0.
\end{equation}
One can show that neither $z(u)\equiv R_1$ nor
$z(u)\equiv R_2$ is
a solution of (\ref{27}) for any $t\in [0,1]$.
Therefore, by the strong maximum principle,
any solution
$z\in C^{4,\a}_a({\Bbb S}^n)$ of (\ref{27})
satisfying $R_1\le z(u)\le R_2$ for all $u\in {\Bbb S}^n$ must
satisfy the strict 
inequalities
\begin{equation}\label{28}
R_1<z(u)<R_2 \quad\mbox{for all~} u\in {\Bbb S}^n.
\end{equation}
By using the $C^1$-estimates in \cite{BLO}, Theorem
\ref{t2}, the
result of Evans and Krylov, and Schauder theory
for second
order uniformly elliptic equations one can obtain
\begin{equation}\label{29}
\|z\|_{C^{4,\a}({\Bbb S}^n)}<C
\end{equation}
for any solution $z\in C^{4,\a}_a({\Bbb S}^n)$ of
(\ref{27}) satisfying (\ref{28}),
where $C$ is a constant depending only on $k$, $n$,
$R_1$, $R_2$
and $\|\psi\|_{C^{2,\a}(\overline{\Omega})}$.

We can choose a constant $\delta>0$ depending on $k$,
$n$, $R_1$,
$R_2$ and $C$ such that
$$
\delta\le \opsi^t(u, z(u))\le \delta^{-1} \quad
\mbox{for~} u\in
{\Bbb S}^n,
$$
where $0\le t\le 1$ and $z\in C^{4,\a}({\Bbb S}^n)$
satisfying
(\ref{28}) and (\ref{29}). Consequently we can find an
open set
$V$ of $\Ga_k$ satisfying $\overline{V}\subset \Ga_k$
such that
$\lambda({\bf B})\in V$ for any $z\in C^{4,\a}_a({\Bbb
S}^n)$
satisfying (\ref{28}), (\ref{29}) and $\delta\le
F({\bf B})\le
\delta^{-1}$. Now we define an open bounded subset
$O^*$ of
$C^{4,\a}({\Bbb S}^n)$ by
$$
O^*:=\{z\in C^{4,\a}({\Bbb S}^n): z \mbox{~satisfies
(\ref{28}),
(\ref{29}) and ~} \lambda({\bf B})\in V\}
$$
One can show that
$$
\Psi(\cdot,t)^{-1}(0)\cap \p O^*=\emptyset \quad
\mbox{for~} 0\le
t\le 1
$$
when $\Psi(\cdot, t)$ are viewed as maps from
$\overline{O}^*\subset C^{4,\a}({\Bbb S}^n)$ to
$C^{2,\a}({\Bbb
S}^n)$. Therefore, the degree $\deg(\Psi(\cdot,t),
O^*, 0)$ is
defined for all $0\le t\le 1$ and is independent of
$t$; see \cite{L}.

 Comparing a solution with
 spheres $z\equiv constant$ and using the maximum principle as usual, 
we know  that $z_0(u)\equiv \overline{R}$
is the unique  solution in $O^*$  of the equation
$\Psi(z,0)=0$. 
Clearly, the linearized operator
$\Psi_z(z_0,0)$
 is of the form
$$
\Psi_z(z_0,0)=-a^{ij}(u)\nabla_{ij}+b^i(u)\nabla_i+c(u),
$$
where $(a^{ij}(u))$ is positive definite.
Since
$$
\Psi(sz_0, 0)=\coth(sz_0)-A^\epsilon \coth^{1+\epsilon}(sz_0),
$$
we have, in view of $A\coth(\overline R)=1$,
\begin{eqnarray*}
\overline  R c(u)&=& \Psi_z(z_0,0)(z_0)=
\frac{d}{ds}\bigg|_{ s=1} \Psi(sz_0, 0)=
\frac{d}{d\rho}\bigg|_{ \rho=\overline  R}
[\coth(\rho)-A^\epsilon \coth^{1+\epsilon}(\rho)]
\\
&=& -\epsilon \frac{d}{d\rho}\bigg|_{ \rho=\overline  R}
\coth(\rho)>0.
\end{eqnarray*}
Thus $\Psi_z(z_0,0)$ is an invertible operator from
$C^{4,\a}({\Bbb
S}^n)$ to $C^{2,\a}({\Bbb S}^n)$.
It follows, as in 
  \cite{LO},
that
$$
\deg(\Psi(\cdot, 1), O^*, 0)=\deg(\Psi(\cdot,0), O^*,
0)\ne 0.
$$
Therefore, the equation
$$
\Psi(z,1)=0,\quad z\in O^*
$$
has at least one solution. This completes the proof of
Theorem
\ref{t1}.

\end{document}